\documentclass{amsart}
\usepackage{amsfonts}

\usepackage{graphicx,hyperref}
\usepackage{amscd}

\newtheorem{theorem}{Theorem}
\theoremstyle{plain}

\newtheorem{corollary}{Corollary}

\newtheorem{remark}{Remark}

\numberwithin{equation}{section}

\input{tcilatex}

\begin{document}
\title[Median Principle for Inequalities]
{The Median Principle for Inequalities and Applications}
\author{S.S. Dragomir}
\address{School of Communications and Informatics\\
Victoria University of Technology\\
PO Box 14428, MCMC 8001\\
Victoria, Australia.}
\email{sever@matilda.vu.edu.au}
\urladdr{http://rgmia.vu.edu.au/SSDragomirWeb.html}
\date{May 23, 2002.}

\begin{abstract}
The ``Median Principle'' for different integral inequalities of Gr\"{u}ss
and Ostrowski type is applied.
\end{abstract}

\maketitle

\section{Introduction}

There are many mathematical inequalities whose right hand side may be
expressed in terms of the sup-norm of a certain derivative for the involved
functions.

For instance, in Numerical Analysis, the integral of a function $f:\left[ a,b%
\right] \rightarrow \mathbb{R}$ may be represented by
\begin{equation}
\int_{a}^{b}f\left( t\right) dt=A_{n}\left( I_{n},f\right) +R_{n}\left(
I_{n},f^{\left( r\right) }\right) ,  \label{1.1}
\end{equation}
where $A_{n}\left( I_{n},f\right) $ is the \textit{quadrature rule} defined
on a given division
\begin{equation*}
I_{n}:a=x_{0}<x_{1}<\cdots <x_{n-1}<x_{n}=b,
\end{equation*}
of the interval $\left[ a,b\right] $ and $R_{n}\left( I_{n},f^{\left(
r\right) }\right) $ is the \textit{remainder}, usually expressed in the
integral form
\begin{equation}
R_{n}\left( I_{n},f^{\left( r\right) }\right) =\int_{a}^{b}K_{n}\left(
I_{n},t\right) f^{\left( r\right) }\left( t\right) dt,  \label{1.2}
\end{equation}
where $K_{n}\left( I_{n},\cdot \right) :\left[ a,b\right] \rightarrow
\mathbb{R}$ is an appropriate \textit{Peano kernel }and $f^{\left( r\right)
} $ is the $r-$th derivative of $f$ assumed to be essentially bounded on $%
\left[ a,b\right] .$

If the integral
\begin{equation*}
\int_{a}^{b}\left| K_{n}\left( I_{n},t\right) \right| dt
\end{equation*}
can be exactly computed or bounded above by different techniques, then we
have the \textit{error estimate}
\begin{equation}
\left| R_{n}\left( I_{n},f^{\left( r\right) }\right) \right| \leq \left\|
f^{\left( r\right) }\right\| _{\left[ a,b\right] ,\infty }\int_{a}^{b}\left|
K_{n}\left( I_{n},t\right) \right| dt,  \label{1.3}
\end{equation}
where
\begin{equation*}
\left\| f^{\left( r\right) }\right\| _{\left[ a,b\right] ,\infty
}:=ess\sup\limits_{t\in \left[ a,b\right] }\left| f^{\left( r\right) }\left(
t\right) \right| ,
\end{equation*}
that provides a large class of examples of inequalities mentioned above.

In Analytic Inequalities Theory, the results such as Ostrowski's inequality
\begin{multline}
\left| f\left( x\right) -\frac{1}{b-a}\int_{a}^{b}f\left( t\right) dt\right|
\label{1.4} \\
\leq \left[ \frac{1}{4}+\left( \frac{x-\frac{a+b}{2}}{b-a}\right) ^{2}\right]
\left( b-a\right) \left\| f^{\prime }\right\| _{\left[ a,b\right] ,\infty
},\;\;x\in \left[ a,b\right] ,
\end{multline}
provided $f$ is absolutely continuous with $f^{\prime }\in L_{\infty }\left[
a,b\right] ,$ or \v{C}eby\v{s}ev's inequality
\begin{multline}
\left| \frac{1}{b-a}\int_{a}^{b}f\left( t\right) g\left( t\right) dt-\frac{1%
}{b-a}\int_{a}^{b}f\left( t\right) dt\cdot \frac{1}{b-a}\int_{a}^{b}g\left(
t\right) dt\right|   \label{1.5} \\
\leq \frac{1}{12}\left( b-a\right) ^{2}\left\| f^{\prime }\right\| _{\left[
a,b\right] ,\infty }\left\| g^{\prime }\right\| _{\left[ a,b\right] ,\infty
},
\end{multline}
provided $f$ and $g$ are absolutely continuous with $f^{\prime },g^{\prime
}\in L_{\infty }\left[ a,b\right] ,$ are other natural examples.

Since, in order to estimate $\left\| f^{\left( r\right) }\right\| _{\left[
a,b\right] ,\infty }$, in practice it is usually necessary to find the
quantities
\begin{equation*}
M_{r}:=\sup\limits_{t\in \left[ a,b\right] }f^{\left( r\right) }\left(
t\right) \text{ \hspace{0.05in}and \hspace{0.05in}}m_{r}:=\inf\limits_{t\in %
\left[ a,b\right] }f^{\left( r\right) }\left( t\right)
\end{equation*}
(as $\left\| f^{\left( r\right) }\right\| _{\left[ a,b\right] ,\infty }=\max
\left\{ \left| M_{r}\right| ,\left| m_{r}\right| \right\} $), the knowledge
of $\left\| f^{\left( r\right) }\right\| _{\left[ a,b\right] ,\infty }$ may
be as difficult as the knowledge of $M_{r}$ and $m_{r}.$

Consequently, it is a natural problem in trying to establish inequalities
where instead of $\left\| f^{\left( r\right) }\right\| _{\left[ a,b\right]
,\infty },$ one would have the positive quantity $M_{r}-m_{r}.$ We also must
note that for functions whose derivatives $f^{\left( r\right) }$ have a
``modest variation'', the quantity $M_{r}-m_{r}$ may usually be a lot
smaller than $\left\| f^{\left( r\right) }\right\| _{\left[ a,b\right]
,\infty }.$

We note that there are many examples of inequalities where the bounds are
expressed in terms of $M_{r}-m_{r}$, from which we would like to mention
only the celebrated result due to Gr\"{u}ss
\begin{multline}
\left| \frac{1}{b-a}\int_{a}^{b}f\left( t\right) g\left( t\right) dt-\frac{1%
}{b-a}\int_{a}^{b}f\left( t\right) dt\cdot \frac{1}{b-a}\int_{a}^{b}g\left(
t\right) dt\right|   \label{1.6} \\
\leq \frac{1}{4}\left( M-m\right) \left( N-n\right) ,
\end{multline}
provided $f,g\in L\left[ a,b\right] $ with
\begin{equation}
m\leq f\leq M,\;\;n\leq g\leq N\text{ \hspace{0.05in}a.e. on }\left[ a,b%
\right] .  \label{1.7}
\end{equation}

It is the main purpose of this paper to point out a general strategy for
transforming an inequality whose right hand side is expressed in terms of $%
\left\| f^{\left( r\right) }\right\| _{\left[ a,b\right] ,\infty }$ into an
inequality for which the same side will be expressed in terms of the
quantity $M_{r}-m_{r}>0.$ We call this method the ``Median Principle''. A
formal statement of this method is provided in the next section.
Applications for some well-known inequalities are given as well.

\section{The Median Principle}

Consider the class of polynomials
\begin{equation*}
\mathcal{P}_{n}^{0}:=\left\{ P_{n}|P_{n}\left( x\right)
=x^{n}+a_{1}x^{n-1}+a_{2}x^{n-2}+\cdots +a_{n-1}x+a_{n},\;a_{i}\in \mathbb{R}%
\right\} .
\end{equation*}
The following result, that will be called the ``\textit{Median Principle'', }%
holds.

\begin{theorem}
\label{t2.1}Let $f:\left[ a,b\right] \subset \mathbb{R}\rightarrow \mathbb{R}
$ be a function so that $f^{\left( n-1\right) }$ is absolutely continuous
and $f^{\left( n\right) }\in L_{\infty }\left[ a,b\right] .$ Assume that the
following inequality holds
\begin{equation}
L\left( f,f^{\left( 1\right) },\dots ,f^{\left( n-1\right) },f^{\left(
n\right) };a,b\right) \leq R\left( \left\| f^{\left( n\right) }\right\| _{
\left[ a,b\right] ,\infty };a,b\right) ,  \label{2.1}
\end{equation}
where $L\left( \cdot ,\cdot ,\cdots ;a,b\right) :\mathbb{R}^{\left(
n+1\right) }\rightarrow \mathbb{R}$ is a general function, $R:[0,\infty
)\rightarrow \mathbb{R}$ and $R$ is monotonic nondecreasing on $[0,\infty ).$

If $g:\left[ a,b\right] \rightarrow \mathbb{R}$ is such that $g^{\left(
n-1\right) }$ is absolutely continuous and
\begin{equation}
-\infty <\gamma \leq g^{\left( n-1\right) }\left( x\right) \leq \Gamma
<\infty \text{ for a.e. }x\in \left[ a,b\right] ,  \label{2.2}
\end{equation}
then one has the inequality
\begin{multline}
\sup\limits_{P_{n}\in \mathcal{P}_{n}^{0}}L\left( g-\frac{\gamma +\Gamma }{2}%
P_{n},g^{\left( 1\right) }-\frac{\gamma +\Gamma }{2}P_{n}^{\left( 1\right)
},\right.   \label{2.3.a} \\
\dots ,\left. g^{\left( n-1\right) }-\frac{\gamma +\Gamma }{2}P_{n}^{\left(
n-1\right) },g^{\left( n\right) }-\frac{\gamma +\Gamma }{2};a,b\right)  \\
\leq R\left( \frac{\Gamma -\gamma }{2};a,b\right) .
\end{multline}
\end{theorem}

\begin{proof}
Let $P_{n}\in \mathcal{P}_{n}^{\left( 0\right) }$ and define $f:\left[ a,b%
\right] \rightarrow \mathbb{R}$, $f\left( x\right) =g\left( x\right) -\frac{%
\gamma +\Gamma }{2}P_{n}^{\left( n\right) }\left( x\right) ,$ where $g$
satisfies (\ref{2.2}).

Obviously, $f^{\left( n\right) }\in L_{\infty }\left[ a,b\right] $ and
\begin{equation*}
f^{\left( n\right) }\left( x\right) =g^{\left( n\right) }\left( x\right) -%
\frac{\gamma +\Gamma }{2}P_{n}^{\left( n\right) }\left( x\right) =g^{\left(
n\right) }\left( x\right) -\frac{\gamma +\Gamma }{2}.
\end{equation*}
Also,
\begin{equation*}
\left| f^{\left( n\right) }\left( x\right) \right| =\left| g^{\left(
n\right) }\left( x\right) -\frac{\gamma +\Gamma }{2}\right| \leq \frac{%
\Gamma -\gamma }{2}
\end{equation*}
giving
\begin{equation*}
\left\| f^{\left( n\right) }\right\| _{\left[ a,b\right] ,\infty }\leq \frac{%
\Gamma -\gamma }{2}.
\end{equation*}
Using the inequality (\ref{2.1}) and the monotonicity of $R,$ we deduce the
desired inequality (\ref{2.3.a}).
\end{proof}

\begin{remark}
\label{r2.1}Similar results may be obtained if the function $L$ or/and $R$
depend on other function $h,l,$ etc.
\end{remark}

The example provided in the next sections will show how the ``Median
Principle'' works in practice.

\section{Inequalities of the $0^{\text{th}}-$Degree}

An inequality of the form
\begin{equation}
L\left( h;a,b\right) \leq R\left( \left\| h\right\| _{\left[ a,b\right]
,\infty };a,b\right) ,  \label{3.1}
\end{equation}
i.e., no derivatives of the function $f$ are involved, is said to be of the $%
0^{th}$ degree.

For example, the following inequality:
\begin{equation}
\left| \int_{a}^{b}h\left( x\right) l\left( x\right) dx\right| \leq \left\|
h\right\| _{\left[ a,b\right] ,\infty }\int_{a}^{b}\left| l\left( x\right)
\right| dx,  \label{3.1a}
\end{equation}
provided $h\in L_{\infty }\left[ a,b\right] $ and $l\in L_{1}\left[ a,b%
\right] $, and
\begin{equation}
\left| \int_{a}^{b}h\left( x\right) du\left( x\right) \right| \leq \left\|
h\right\| _{\left[ a,b\right] ,\infty }\bigvee_{a}^{b}\left( u\right) ,
\label{3.2}
\end{equation}
provided $h\in C\left[ a,b\right] $ (the class of continuous functions) and $%
u\in BV\left[ a,b\right] $ (the class of functions of bounded variation),
are inequalities of $0^{th}-$degree. As a generalisation of (\ref{3.2}), if $%
h,l\in C\left[ a,b\right] $ and $u\in BV\left[ a,b\right] ,$ then also
\begin{equation}
\left| \int_{a}^{b}h\left( x\right) l\left( x\right) du\left( x\right)
\right| \leq \left\| h\right\| _{\left[ a,b\right] ,\infty }\left\|
l\right\| _{\left[ a,b\right] ,\infty }\bigvee_{a}^{b}\left( u\right) .
\label{3.3}
\end{equation}
Here and in (\ref{3.2}), $\bigvee_{a}^{b}\left( u\right) $ denotes the total
variation of $u$ in $\left[ a,b\right] .$

The following result holds.

\begin{theorem}
\label{t3.1}Let $f,l:\left[ a,b\right] \rightarrow \mathbb{R}$ be such that
there exists the constants $m,M\in \mathbb{R}$ with the property
\begin{equation}
-\infty <m\leq f\left( x\right) \leq M<\infty \text{ for a.e. }x\in \left[
a,b\right] ,  \label{3.4}
\end{equation}
and $l\in L_{1}\left[ a,b\right] ,$ such that
\begin{equation}
\int_{a}^{b}l\left( x\right) dx=0.  \label{3.5}
\end{equation}
Then we have the inequality
\begin{equation}
\left| \int_{a}^{b}f\left( x\right) l\left( x\right) dx\right| \leq \frac{1}{%
2}\left( M-m\right) \int_{a}^{b}\left| l\left( x\right) \right| dx.
\label{3.6}
\end{equation}
The constant $\frac{1}{2}$ is sharp.
\end{theorem}

\begin{proof}
Using the ``\textit{Median Principle'' }for the inequality (\ref{3.1a}) we
have
\begin{equation}
\left| \int_{a}^{b}\left( f\left( x\right) -\frac{m+M}{2}\right) l\left(
x\right) dx\right| \leq \frac{1}{2}\left( M-m\right) \int_{a}^{b}\left|
l\left( x\right) \right| dx  \label{3.7}
\end{equation}
and since $\int_{a}^{b}l\left( x\right) dx=0,$ we deduce (\ref{3.6}).

Now, assume that (\ref{3.6}) holds with a constant $C>0,$ i.e.,
\begin{equation}
\left| \int_{a}^{b}f\left( x\right) l\left( x\right) dx\right| \leq C\left(
M-m\right) \int_{a}^{b}\left| l\left( x\right) \right| dx.  \label{3.8}
\end{equation}
If we choose $f=l$ and $f:\left[ a,b\right] \rightarrow \mathbb{R}$ where
\begin{equation*}
f\left( x\right) =\left\{
\begin{array}{ll}
-1, & x\in \left[ a,\frac{a+b}{2}\right] \\
&  \\
1, & x\in \left( \frac{a+b}{2},b\right] .
\end{array}
\right.
\end{equation*}
Then
\begin{eqnarray*}
&&\left. \left| \int_{a}^{b}f\left( x\right) l\left( x\right) dx\right|
=b-a,\right. \\
&&\left. \int_{a}^{b}\left| l\left( x\right) \right| dx=b-a,\right. \\
&&\left. m=-1,\;M=1\right.
\end{eqnarray*}
and thus, by (\ref{3.8}) we deduce $C\geq \frac{1}{2},$ and the theorem is
then proved.
\end{proof}

\begin{corollary}
\label{c3.2}Let $f,g:\left[ a,b\right] \rightarrow \mathbb{R}$ be such that $%
f$ satisfies (\ref{3.4}) and $g\in L_{1}\left[ a,b\right] .$ Then we have
the inequalities
\begin{multline}
\left| \frac{1}{b-a}\int_{a}^{b}f\left( x\right) g\left( x\right) dx-\frac{1%
}{b-a}\int_{a}^{b}f\left( x\right) dx\cdot \frac{1}{b-a}\int_{a}^{b}g\left(
x\right) dx\right|  \label{3.9} \\
\leq \frac{1}{2}\left( M-m\right) \frac{1}{b-a}\int_{a}^{b}\left| g\left(
x\right) -\frac{1}{b-a}\int_{a}^{b}g\left( y\right) dy\right| dx.
\end{multline}
The constant $\frac{1}{2}$ is sharp in (\ref{3.9}).
\end{corollary}

\begin{proof}
Follows by (\ref{3.6}) on choosing $l\left( x\right) =g\left( x\right) -%
\frac{1}{b-a}\int_{a}^{b}g\left( y\right) dy.$
\end{proof}

\begin{remark}
\label{r3.3}The inequality (\ref{3.5}) was proved in a different, more
complicated, way in \cite{CS}. Generalisations for abstract Lebesgue
integrals, the weighted and discrete cases were obtained in \cite{CD}.
\end{remark}

The following result also holds.

\begin{theorem}
\label{t3.4}Let $f:\left[ a,b\right] \rightarrow \mathbb{R}$ be a continuous
function with the property that
\begin{equation}
-\infty <m\leq f\left( x\right) \leq M<\infty \text{ for a.e. }x\in \left[
a,b\right] ,  \label{3.10}
\end{equation}
and $u\in BV\left[ a,b\right] $ with the property that
\begin{equation}
u\left( a\right) =u\left( b\right) .  \label{3.11}
\end{equation}
Then we have the inequality:
\begin{equation}
\left| \int_{a}^{b}f\left( x\right) du\left( x\right) \right| \leq \frac{1}{2%
}\left( M-m\right) \bigvee_{a}^{b}\left( u\right) .  \label{3.12}
\end{equation}
The constant $\frac{1}{2}$ is sharp.
\end{theorem}

\begin{proof}
Using the ``Median Principle'' for the inequality (\ref{3.2}) we have
\begin{equation}
\left| \int_{a}^{b}\left( f\left( x\right) -\frac{m+M}{2}\right) du\left(
x\right) \right| \leq \frac{1}{2}\left( M-m\right) \bigvee_{a}^{b}\left(
u\right) .  \label{3.13}
\end{equation}
Since $\int_{a}^{b}du\left( x\right) =u\left( b\right) -u\left( a\right) =0,$
from (\ref{3.13}) we deduce (\ref{3.12}).

Now, assume that the inequality (\ref{3.12}) holds with a constant $D>0.$
That is,
\begin{equation}
\left| \int_{a}^{b}f\left( x\right) du\left( x\right) \right| \leq D\left(
M-m\right) \bigvee_{a}^{b}\left( u\right) .  \label{3.14}
\end{equation}
Consider $a=0,$ $b=1,$ $f:\left[ 0,1\right] \rightarrow \mathbb{R}$, $%
f\left( x\right) =x,$ and $u:\left[ 0,1\right] \rightarrow \mathbb{R}$ given
by
\begin{equation*}
u\left( x\right) =\left\{
\begin{array}{ll}
0 & \text{if \hspace{0.05in}}x=0\text{ or }x=1. \\
&  \\
1 & \text{if \hspace{0.05in}}x\in \left( 0,1\right) .
\end{array}
\right.
\end{equation*}
We have
\begin{equation*}
\int_{0}^{1}f\left( x\right) du\left( x\right) =u\left( x\right) f\left(
x\right) \big|_{0}^{1}-\int_{0}^{1}u\left( x\right) df\left( x\right)
=-\int_{0}^{1}u\left( x\right) dx=-1,
\end{equation*}
\begin{equation*}
M=1,\;m=0
\end{equation*}
and
\begin{equation*}
\bigvee_{a}^{b}\left( u\right) =2.
\end{equation*}
Then, by (\ref{3.14}) we deduce $2D\geq 1$ giving $D\geq \frac{1}{2},$ and
the theorem is thus proved.
\end{proof}

Another result generalizing the above ones also holds.

\begin{theorem}
\label{t3.5}Let $f,l:\left[ a,b\right] \rightarrow \mathbb{R}$ be continuous
and $f$ is such that the condition (\ref{3.10}) holds. If $u\in BV\left( %
\left[ a,b\right] \right) $ is such that
\begin{equation}
\int_{a}^{b}l\left( x\right) du\left( x\right) =0,  \label{3.15}
\end{equation}
then we have the inequality:
\begin{equation}
\left| \int_{a}^{b}f\left( x\right) l\left( x\right) du\left( x\right)
\right| \leq \frac{1}{2}\left( M-m\right) \left\| l\right\| _{\left[ a,b%
\right] ,\infty }\bigvee_{a}^{b}\left( u\right) .  \label{3.16}
\end{equation}
The constant $\frac{1}{2}$ in (\ref{3.16}) is sharp.
\end{theorem}

\begin{proof}
Follows by the ``Median Principle'' applied for the inequality (\ref{3.3}).
The sharpness of the constant follows by Theorem \ref{t3.4} on choosing $%
l=1. $
\end{proof}

As a corollary of the above result, we may state the following Gr\"{u}ss
type inequality.

\begin{corollary}
\label{c3.6}Let $f,g\in C\left[ a,b\right] $ and $f$ is such that (\ref{3.10}%
) holds. If $u\in BV\left[ a,b\right] $ and $u\left( b\right) \neq u\left(
a\right) ,$ then one has the inequality
\begin{multline}
\left| \frac{1}{u\left( b\right) -u\left( a\right) }\int_{a}^{b}f\left(
x\right) g\left( x\right) du\left( x\right) \right.   \label{3.17.a} \\
-\left. \frac{1}{u\left( b\right) -u\left( a\right) }\int_{a}^{b}f\left(
x\right) du\left( x\right) \cdot \frac{1}{u\left( b\right) -u\left( a\right)
}\int_{a}^{b}g\left( x\right) du\left( x\right) \right|  \\
\leq \frac{1}{2}\left( M-m\right) \frac{1}{\left| u\left( b\right) -u\left(
a\right) \right| }\left\| g-\frac{1}{u\left( b\right) -u\left( a\right) }%
\int_{a}^{b}g\left( y\right) du\left( y\right) \right\| _{\left[ a,b\right]
,\infty }\bigvee_{a}^{b}\left( u\right) .
\end{multline}
The constant $\frac{1}{2}$ is sharp in (\ref{3.17.a}).
\end{corollary}

\begin{proof}
We choose in Theorem \ref{t3.5}, $l:\left[ a,b\right] \rightarrow \mathbb{R}%
, $%
\begin{equation*}
l\left( x\right) =g\left( x\right) -\frac{1}{u\left( b\right) -u\left(
a\right) }\int_{a}^{b}g\left( y\right) du\left( y\right) ,\;\;x\in \left[ a,b%
\right] .
\end{equation*}
Then, obviously
\begin{equation*}
\int_{a}^{b}l\left( x\right) du\left( x\right) =0,
\end{equation*}
and by (\ref{3.16}) we deduce (\ref{3.17.a}).

To prove the sharpness of the constant $\frac{1}{2}$ in (\ref{3.17.a}), we
assume that it holds with a constant $C>0,$ i.e.,
\begin{multline}
\left| \frac{1}{u\left( b\right) -u\left( a\right) }\int_{a}^{b}f\left(
x\right) g\left( x\right) du\left( x\right) \right.   \label{3.18} \\
-\left. \frac{1}{u\left( b\right) -u\left( a\right) }\int_{a}^{b}f\left(
x\right) du\left( x\right) \cdot \frac{1}{u\left( b\right) -u\left( a\right)
}\int_{a}^{b}g\left( x\right) du\left( x\right) \right|  \\
\leq C\left( M-m\right) \frac{1}{\left| u\left( b\right) -u\left( a\right)
\right| }\left\| g-\frac{1}{u\left( b\right) -u\left( a\right) }%
\int_{a}^{b}g\left( y\right) du\left( y\right) \right\| _{\left[ a,b\right]
,\infty }\bigvee_{a}^{b}\left( u\right) .
\end{multline}
Let us choose $f=g,$ $f:\left[ a,b\right] \rightarrow \mathbb{R}$, with $%
f\left( t\right) =t$ and $u:\left[ a,b\right] \rightarrow \mathbb{R}$ given
by
\begin{equation*}
u\left( t\right) =\left\{
\begin{array}{ll}
-1, & \text{if }t=a \\
&  \\
0 & \text{if }t\in \left( a,b\right) , \\
&  \\
1, & \text{if }t=b.
\end{array}
\right.
\end{equation*}
Then
\begin{eqnarray*}
\frac{1}{u\left( b\right) -u\left( a\right) }\int_{a}^{b}f\left( x\right)
g\left( x\right) du\left( x\right)  &=&\frac{1}{2}\int_{a}^{b}u^{2}du\left(
x\right)  \\
&=&\frac{1}{2}\left[ x^{2}u\left( x\right) \bigg|_{a}^{b}-2\int_{a}^{b}xu%
\left( x\right) dx\right]  \\
&=&\frac{b^{2}+a^{2}}{2},
\end{eqnarray*}
\begin{eqnarray*}
\int_{a}^{b}f\left( x\right) du\left( x\right)  &=&\int_{a}^{b}g\left(
x\right) du\left( x\right) =\frac{1}{2}\int_{a}^{b}xdu\left( x\right)  \\
&=&xu\left( x\right) \bigg|_{a}^{b}-\int_{a}^{b}u\left( x\right) dx=b+a,
\end{eqnarray*}
\begin{equation*}
\left\| g-\frac{1}{u\left( b\right) -u\left( a\right) }\int_{a}^{b}g\left(
s\right) ds\right\| _{\infty }=\sup\limits_{x\in \left[ a,b\right] }\left| x-%
\frac{a+b}{2}\right| =\frac{b-a}{2},
\end{equation*}
\begin{equation*}
\bigvee_{a}^{b}\left( u\right) =2,\;M=b,\;\;m=a
\end{equation*}
and thus, by (\ref{3.18}), we get
\begin{equation*}
\left| \frac{a^{2}+b^{2}}{2}-\frac{\left( a+b\right) ^{2}}{4}\right| \leq
C\left( b-a\right) \frac{1}{2}\cdot \frac{\left( b-a\right) }{2}\cdot 2,
\end{equation*}
giving $C\geq \frac{1}{2},$ and the corollary is proved.
\end{proof}

For other results of this type see \cite{SSD2}.

\section{Inequalities of the $1^{\text{st}}-$Degree}

An inequality that contains at most the first derivative of the involved
functions will be called an inequality of the $1^{st}-$degree.

For example, Ostrowski's inequality
\begin{multline}
\left| h\left( x\right) -\frac{1}{b-a}\int_{a}^{b}h\left( t\right) dt\right|
\label{4.1} \\
\leq \left[ \frac{1}{4}+\left( \frac{x-\frac{a+b}{2}}{b-a}\right) ^{2}\right]
\left\| h^{\prime }\right\| _{\left[ a,b\right] ,\infty }\left( b-a\right)
,\;\;x\in \left[ a,b\right] ;
\end{multline}
provided $h$ is absolutely continuous and $h^{\prime }\in L_{\infty }\left[
a,b\right] ,$ is such an inequality,

Also, the generalised-trapezoid inequality:
\begin{multline}
\left| \frac{\left( x-a\right) h\left( a\right) +\left( b-x\right) h\left(
b\right) }{b-a}-\frac{1}{b-a}\int_{a}^{b}h\left( t\right) dt\right|
\label{4.2} \\
\leq \left[ \frac{1}{4}+\left( \frac{x-\frac{a+b}{2}}{b-a}\right) ^{2}\right]
\left\| h^{\prime }\right\| _{\left[ a,b\right] ,\infty }\left( b-a\right)
,\;\;x\in \left[ a,b\right] ;
\end{multline}
provided $h$ is absolutely continuous, $h^{\prime }\in L_{\infty }\left[ a,b%
\right] ,$ is another example of such an inequality.

In both the inequalities above, the constant $\frac{1}{4}$ is sharp in the
sense that it cannot be replaced by a smaller constant.

If one would like examples of such inequalities for two functions, the
following Ostrowski's inequality obtained in \cite{AO} is the most suitable
\begin{multline}
\left| \frac{1}{b-a}\int_{a}^{b}h\left( x\right) l\left( x\right) dx-\frac{1%
}{b-a}\int_{a}^{b}h\left( x\right) dx\cdot \frac{1}{b-a}\int_{a}^{b}l\left(
x\right) dx\right|  \label{4.3} \\
\leq \frac{1}{8}\left( b-a\right) \left( M-m\right) \left\| h^{\prime
}\right\| _{\left[ a,b\right] ,\infty },
\end{multline}
provided $-\infty <m\leq h\left( x\right) \leq M<\infty $ for a.e. $x\in %
\left[ a,b\right] ,$ and $l$ is absolutely continuous and such that $%
l^{\prime }\in L_{\infty }\left[ a,b\right] .$ The constant $\frac{1}{8}$ is
sharp.

Another example of such an inequality is the \v{C}eby\v{s}ev one
\begin{multline}
\left| \frac{1}{b-a}\int_{a}^{b}h\left( x\right) l\left( x\right) dx-\frac{1%
}{b-a}\int_{a}^{b}h\left( x\right) dx\cdot \frac{1}{b-a}\int_{a}^{b}l\left(
x\right) dx\right|  \label{4.4} \\
\leq \frac{1}{12}\left( b-a\right) ^{2}\left\| h^{\prime }\right\| _{\left[
a,b\right] ,\infty }\left\| l^{\prime }\right\| _{\left[ a,b\right] ,\infty
},
\end{multline}
provided $h,l$ are absolutely continuous and $h^{\prime },l^{\prime }\in
L_{\infty }\left[ a,b\right] .$ The constant $\frac{1}{12}$ here is sharp.

The following perturbed version of Ostrowski's inequality holds.

\begin{theorem}
\label{t4.1}Let $f:\left[ a,b\right] \rightarrow \mathbb{R}$ be an
absolutely continuous function on $\left[ a,b\right] $ such that
\begin{equation}
-\infty <\gamma \leq f^{\prime }\left( x\right) \leq \Gamma <\infty \text{
for a.e. }x\in \left[ a,b\right] .  \label{4.5}
\end{equation}
Then one has the inequality:
\begin{multline}
\left| f\left( x\right) -\frac{1}{b-a}\int_{a}^{b}f\left( t\right) dt-\frac{%
\gamma +\Gamma }{2}\left( x-\frac{a+b}{2}\right) \right|  \label{4.6} \\
\leq \frac{1}{2}\left[ \frac{1}{4}+\left( \frac{x-\frac{a+b}{2}}{b-a}\right)
^{2}\right] \left( \Gamma -\gamma \right) \left( b-a\right)
\end{multline}
for any $x\in \left[ a,b\right] .$

The constant $\frac{1}{2}$ is sharp.
\end{theorem}

\begin{proof}
Consider the function $h\left( x\right) =f\left( x\right) -\frac{\gamma
+\Gamma }{2}x,$ $x\in \left[ a,b\right] .$ Applying the ``Median Principle''
for Ostrowski's inequality, we get
\begin{multline*}
\left| h\left( x\right) -\frac{\gamma +\Gamma }{2}x-\frac{1}{b-a}%
\int_{a}^{b}\left( h\left( t\right) -\frac{\gamma +\Gamma }{2}t\right)
dt\right| \\
\leq \frac{1}{2}\left( \Gamma -\gamma \right) \left[ \frac{1}{4}+\left(
\frac{x-\frac{a+b}{2}}{b-a}\right) ^{2}\right] \left( b-a\right) ,
\end{multline*}
which is clearly equivalent to (\ref{4.6}).

The sharpness of the constant follows by the sharpness of Ostrowski's
inequality on choosing $\gamma =-\left\| f^{\prime }\right\| _{\left[ a,b%
\right] ,\infty },\;\Gamma =\left\| f^{\prime }\right\| _{\left[ a,b\right]
,\infty }.$ We omit the details.
\end{proof}

\begin{remark}
For a different proof of this fact, see \cite{SSD1}.
\end{remark}

Now, we may give a perturbed version of the generalised trapezoid inequality
(\ref{4.2}) as well (see also \cite{SSD1}).

\begin{theorem}
\label{t4.2}Let $f$ be as in Theorem \ref{t4.1}. Then one has the inequality
\begin{multline*}
\left| \frac{\left( x-a\right) f\left( a\right) +\left( b-x\right) f\left(
b\right) }{b-a}-\frac{\gamma +\Gamma }{2}\left( x-\frac{a+b}{2}\right) -%
\frac{1}{b-a}\int_{a}^{b}h\left( t\right) dt\right| \\
\leq \frac{1}{2}\left( \Gamma -\gamma \right) \left[ \frac{1}{4}+\left(
\frac{x-\frac{a+b}{2}}{b-a}\right) ^{2}\right] \left( b-a\right) ,
\end{multline*}
for any $x\in \left[ a,b\right] .$

The constant $\frac{1}{2}$ is sharp.
\end{theorem}

The proof follows by the inequality (\ref{4.2}) and we omit the details.

Now, we are able to point out the following perturbation of the second
Ostrowski's inequality (\ref{4.3}).

\begin{theorem}
\label{t4.3}Let $f:\left[ a,b\right] \rightarrow \mathbb{R}$ be an
absolutely continuous function on $\left[ a,b\right] $ such that the
derivative $f^{\prime }:\left[ a,b\right] \rightarrow \mathbb{R}$ satisfies
the condition
\begin{equation}
-\infty <\gamma \leq f^{\prime }\left( x\right) \leq \Gamma <\infty \text{
for a.e. }x\in \left[ a,b\right] .  \label{4.7}
\end{equation}
If $g:\left[ a,b\right] \rightarrow \mathbb{R}$ is such that
\begin{equation}
-\infty <m\leq g\left( x\right) \leq M<\infty \text{ for a.e. }x\in \left[
a,b\right] ,  \label{4.8}
\end{equation}
then we have the inequality:
\begin{multline}
\left| \frac{1}{b-a}\int_{a}^{b}f\left( x\right) g\left( x\right) dx-\frac{1%
}{b-a}\int_{a}^{b}f\left( x\right) dx\cdot \frac{1}{b-a}\int_{a}^{b}g\left(
x\right) dx\right.  \label{4.9} \\
\left. -\frac{\gamma +\Gamma }{2}\cdot \frac{1}{b-a}\int_{a}^{b}\left( x-%
\frac{a+b}{2}\right) g\left( x\right) dx\right| \\
\leq \frac{1}{16}\left( b-a\right) \left( M-m\right) \left( \Gamma -\gamma
\right)
\end{multline}
The constant $\frac{1}{16}$ is best possible.
\end{theorem}

\begin{proof}
Consider $h\left( x\right) =f\left( x\right) -\frac{\gamma +\Gamma }{2}x.$
Applying the ``Median Principle'' for the Ostrowski's inequality (\ref{4.3}%
), we have:
\begin{multline}
\left| \frac{1}{b-a}\int_{a}^{b}\left( f\left( x\right) -\frac{\gamma
+\Gamma }{2}x\right) g\left( x\right) dx\right.  \label{4.10} \\
\left. -\frac{1}{b-a}\int_{a}^{b}\left( f\left( x\right) -\frac{\gamma
+\Gamma }{2}x\right) dx\cdot \frac{1}{b-a}\int_{a}^{b}g\left( x\right)
dx\right| \\
\leq \frac{1}{16}\left( b-a\right) \left( M-m\right) \left( \Gamma -\gamma
\right)
\end{multline}
that after some elementary computations is equivalent to (\ref{4.9}).

The sharpness of the constant $\frac{1}{16}$ follows by the fact that the
constant $\frac{1}{8}$ is sharp in (\ref{4.3}) on taking $\gamma =-\left\|
f^{\prime }\right\| _{\left[ a,b\right] ,\infty },\;\Gamma =\left\|
f^{\prime }\right\| _{\left[ a,b\right] ,\infty }.$ We omit the details.
\end{proof}

\section{Inequalities of the $n^{\text{th}}-$Degree}

In \cite{CDR}, the authors proved the following identity
\begin{multline}
\int_{a}^{b}f\left( t\right) dt=\sum_{k=0}^{n-1}\left[ \frac{\left(
b-x\right) ^{k+1}+\left( -1\right) ^{k}\left( x-a\right) ^{k+1}}{\left(
k+1\right) !}\right] f^{\left( k\right) }\left( x\right)  \label{5.1} \\
+\left( -1\right) ^{n}\int_{a}^{b}K_{n}\left( x,t\right) f^{\left( n\right)
}\left( t\right) dt,
\end{multline}
provided $f^{\left( n-1\right) }$ is absolutely continuous on $\left[ a,b%
\right] $ and the kernel $K_{n}:\left[ a,b\right] ^{2}\rightarrow \mathbb{R}$
is given by
\begin{equation}
K_{n}\left( x,t\right) :=\left\{
\begin{array}{ll}
\dfrac{\left( t-a\right) ^{n}}{n!} & \text{if\ \ }a\leq t\leq x\leq b, \\
&  \\
\dfrac{\left( t-b\right) ^{n}}{n!} & \text{if\ \ }a\leq x<t\leq b.
\end{array}
\right.  \label{5.2}
\end{equation}
Using the representation (\ref{5.1}), they proved the following inequality
\begin{multline}
\left| \int_{a}^{b}f\left( t\right) dt-\sum_{k=0}^{n-1}\left[ \frac{\left(
b-x\right) ^{k+1}+\left( -1\right) ^{k}\left( x-a\right) ^{k+1}}{\left(
k+1\right) !}\right] f^{\left( k\right) }\left( x\right) \right|  \label{5.3}
\\
\leq \frac{1}{\left( n+1\right) !}\left\| f^{\left( n\right) }\right\|
_{\infty }\left[ \left( x-a\right) ^{n+1}+\left( b-x\right) ^{n+1}\right]
\end{multline}
for any $x\in \left[ a,b\right] ,$ and in particular, for $x=\frac{a+b}{2}$%
\begin{multline}
\left| \int_{a}^{b}f\left( t\right) dt-\sum_{k=0}^{n-1}\left[ \frac{1+\left(
-1\right) ^{k}}{\left( k+1\right) !}\right] \cdot \frac{\left( b-a\right)
^{k+1}}{2^{k+1}}f^{\left( k\right) }\left( \frac{a+b}{2}\right) \right|
\label{5.4} \\
\leq \frac{1}{2^{n}\left( n+1\right) !}\left\| f^{\left( n\right) }\right\|
_{\infty }\left( b-a\right) ^{n+1}.
\end{multline}

The following result holds.

\begin{theorem}
\label{t5.1}Let $f:\left[ a,b\right] \subset \mathbb{R\rightarrow R}$ be a
function such that the derivative $f^{\left( n-1\right) }$ is absolutely
continuous on $\left[ a,b\right] $ and there exists the constants $\gamma
_{n},\Gamma _{n}\in \mathbb{R}$ so that
\begin{equation}
-\infty <\gamma _{n}\leq f^{\left( n\right) }\left( t\right) \leq \Gamma
_{n}<\infty \text{ for a.e. }t\in \left[ a,b\right] .  \label{5.5}
\end{equation}
Then we have the inequality
\begin{multline}
\left| \int_{a}^{b}f\left( t\right) dt-\sum_{k=0}^{n-1}\left[ \frac{\left(
b-x\right) ^{k+1}+\left( -1\right) ^{k}\left( x-a\right) ^{k+1}}{\left(
k+1\right) !}\right] f^{\left( k\right) }\left( x\right) \right.
\label{5.6} \\
\left. -\left( -1\right) ^{n}\frac{\Gamma _{n}+\gamma _{n}}{2}\left[ \frac{%
\left( x-a\right) ^{n+1}+\left( -1\right) ^{n+1}\left( b-x\right) ^{n+1}}{%
\left( n+1\right) !}\right] \right|  \\
\leq \frac{\Gamma _{n}-\gamma _{n}}{2\left( n+1\right) !}\left[ \left(
x-a\right) ^{n+1}+\left( b-x\right) ^{n+1}\right] .
\end{multline}
\end{theorem}

\begin{proof}
Observe, by (\ref{5.1}), we have that
\begin{align}
\int_{a}^{b}K_{n}\left( x,t\right) dt& =\frac{1}{n!}\left[
\int_{a}^{x}\left( t-a\right) ^{n}dt+\int_{x}^{b}\left( t-b\right) ^{n}dt%
\right]   \label{5.7} \\
& =\frac{\left( x-a\right) ^{n+1}-\left( x-b\right) ^{n+1}}{\left(
n+1\right) !}  \notag \\
& =\frac{\left( x-a\right) ^{n+1}+\left( -1\right) ^{n+1}\left( b-x\right)
^{n+1}}{\left( n+1\right) !}.  \notag
\end{align}
Taking the modulus in (\ref{5.1}) and using the fact that
\begin{equation*}
\left| f^{\left( n+1\right) }\left( t\right) -\frac{\Gamma _{n}+\gamma _{n}}{%
2}\right| \leq \frac{\Gamma _{n}-\gamma _{n}}{2}\text{ \ for a.e. }t\in %
\left[ a,b\right]
\end{equation*}
and
\begin{equation*}
\int_{a}^{b}\left| K_{n}\left( x,t\right) \right| dt=\frac{1}{\left(
n+1\right) !}\left[ \left( x-a\right) ^{n+1}+\left( b-x\right) ^{n+1}\right]
,
\end{equation*}
we easily deduce (\ref{5.6}).
\end{proof}

\begin{corollary}
\label{c5.2}With the assumptions in Theorem \ref{t5.1}, one has the
inequality:
\begin{multline}
\left| \int_{a}^{b}f\left( t\right) dt-\sum_{k=0}^{n-1}\left[ \frac{1+\left(
-1\right) ^{k}}{\left( k+1\right) !}\right] \cdot \frac{\left( b-a\right)
^{k+1}}{2^{k+1}}f^{\left( k\right) }\left( \frac{a+b}{2}\right) \right.
\label{5.8} \\
\left. -\left( -1\right) ^{n}\frac{\Gamma _{n}+\gamma _{n}}{2}\left[ \frac{%
1+\left( -1\right) ^{n+1}}{\left( n+1\right) !}\right] \frac{\left(
b-a\right) ^{n+1}}{2^{n+1}}\right| \\
\leq \frac{\Gamma _{n}-\gamma _{n}}{2^{n+1}\left( n+1\right) !}\left(
b-a\right) ^{n+1}.
\end{multline}
\end{corollary}

In \cite{CDRJ}, the authors also obtained the following identity
\begin{multline}
\int_{a}^{b}f\left( t\right) dt=\sum_{k=0}^{n-1}\frac{1}{\left( k+1\right) !}%
\left[ \left( x-a\right) ^{k+1}f^{\left( k\right) }\left( a\right) +\left(
-1\right) ^{k}\left( b-x\right) ^{k+1}f^{\left( k\right) }\left( b\right) %
\right]  \label{5.9} \\
+\frac{1}{n!}\int_{a}^{b}\left( x-t\right) ^{n}f^{\left( n\right) }\left(
t\right) dt,
\end{multline}
provided $f^{\left( n-1\right) }$ is absolutely continuous on $\left[ a,b%
\right] .$

By the use of this identity, they obtained the inequality
\begin{multline}
\left| \int_{a}^{b}f\left( t\right) dt-\sum_{k=0}^{n-1}\frac{1}{\left(
k+1\right) !}\left[ \left( x-a\right) ^{k+1}f^{\left( k\right) }\left(
a\right) +\left( -1\right) ^{k}\left( b-x\right) ^{k+1}f^{\left( k\right)
}\left( b\right) \right] \right|  \label{5.10} \\
\leq \frac{1}{\left( n+1\right) !}\left\| f^{\left( n\right) }\right\|
_{\infty }\left[ \left( x-a\right) ^{n+1}+\left( b-x\right) ^{n+1}\right] ,
\end{multline}
for any $x\in \left[ a,b\right] .$

In particular, for $x=\frac{a+b}{2}$ we get the inequality:
\begin{multline}
\left| \int_{a}^{b}f\left( t\right) dt-\sum_{k=0}^{n-1}\frac{1}{\left(
k+1\right) !}\left( \frac{b-a}{2}\right) ^{k+1}\left[ f^{\left( k\right)
}\left( a\right) +\left( -1\right) ^{k}f^{\left( k\right) }\left( b\right) %
\right] \right|   \label{5.11} \\
\leq \frac{1}{2^{n}\left( n+1\right) !}\left\| f^{\left( n\right) }\right\|
_{\infty }\left( b-a\right) ^{n+1}.
\end{multline}

Finally, we may state the following result.

\begin{theorem}
\label{t5.3}With the assumptions in Theorem \ref{t5.1}, we have the
inequality
\begin{multline}
\left| \int_{a}^{b}f\left( t\right) dt-\sum_{k=0}^{n-1}\left[ \frac{\left(
b-x\right) ^{k+1}+\left( -1\right) ^{k}\left( x-a\right) ^{k+1}}{\left(
k+1\right) !}\right] f^{\left( k\right) }\left( x\right) \right.
\label{5.12} \\
\left. -\frac{1}{\left( n+1\right) !}\cdot \frac{\Gamma _{n}+\gamma _{n}}{2}%
\left[ \left( x-a\right) ^{n+1}+\left( -1\right) ^{n}\left( b-x\right) ^{n+1}%
\right] \right| \\
\leq \frac{\Gamma _{n}-\gamma _{n}}{2\left( n+1\right) !}\left[ \left(
x-a\right) ^{n+1}+\left( b-x\right) ^{n+1}\right] ,
\end{multline}
for any $x\in \left[ a,b\right] .$
\end{theorem}

In particular, for $x=\frac{a+b}{2},$ we have the corollary

\begin{corollary}
\label{c5.4}With the above assumptions, we have
\begin{multline}
\left| \int_{a}^{b}f\left( t\right) dt-\sum_{k=0}^{n-1}\frac{1}{\left(
k+1\right) !}\left( \frac{b-a}{2}\right) ^{k+1}\left[ f^{\left( k\right)
}\left( a\right) +\left( -1\right) ^{k}f^{\left( k\right) }\left( b\right) %
\right] \right.   \label{5.13} \\
\left. -\frac{\Gamma _{n}+\gamma _{n}}{2}\left[ \frac{1+\left( -1\right) ^{n}%
}{\left( n+1\right) !}\right] \left( \frac{b-a}{2}\right) ^{n+1}\right|  \\
\leq \frac{1}{2^{n+1}\left( n+1\right) !}\left( \Gamma _{n}-\gamma
_{n}\right) \left( b-a\right) ^{n+1}.
\end{multline}
\end{corollary}

The interested reader may find many other examples which can be treated in a
similar fashion. We omit the details.

\end{document}